\documentclass[final,leqno]{siamltex}
\usepackage{epsfig}
\usepackage{amsmath}
\usepackage{amssymb}
\usepackage{tikz}
\usepackage{graphicx}
\usepackage[notcite,notref]{showkeys}
 \usepackage{float}
\newtheorem{algorithm}{Weak Galerkin Algorithm}
\newtheorem{thm1}{Definition}
\newtheorem{remark}{Remark}
\newcommand{\bq}{{\bf q}}
\newcommand{\bn}{{\bf n}}
\newcommand{\bx}{{\bf x}}

\def\T{{\mathcal T}}
\def\E{{\mathcal E}}
\def\Q{{\mathbb Q}}

\def\l{{\langle}}
\def\r{{\rangle}}

\def\bn{{\bf n}}
\def\bq{{\bf q}}

\newcommand{\pT}{{\partial T}}

\def\3bar{{|\hspace{-.02in}|\hspace{-.02in}|}}

  \numberwithin{equation}{section}
\numberwithin{table}{section} \numberwithin{figure}{section}

\title{Numerical investigation on weak Galerkin finite elements}

\author{Xiu Ye\thanks{Department of
Mathematics, University of Arkansas at Little Rock, Little Rock, AR
72204 (xxye@ualr.edu). This research was supported in part by
National Science Foundation Grant DMS-1620016.}
\and
Shangyou Zhang\thanks{Department of
Mathematical Sciences, University of Delaware, Newark, DE 19716 (szhang@udel.edu).}
}

\begin{document}

\maketitle

\begin{abstract}
The weak Galerkin (WG) finite element method is an effective and
flexible general numerical technique for solving partial differential equations.
The novel idea of weak Galerkin finite element methods is on the use of weak functions and
their weak derivatives defined as distributions. Weak functions and weak derivatives
can be approximated by polynomials with various degrees. Different combination of polynomial
spaces generates different weak Galerkin finite elements. The purpose of this paper is
to study stability, convergence and supercloseness of different WG elements
by providing many numerical experiments  recorded in 31 tables. These tables serve two purposes. First it provides
 a detail guide of the performance of different WG elements. Second, the information in the
tables opens new research territory why some WG elements outperform others.
\end{abstract}


\pagestyle{myheadings}

\section{Introduction}\label{Section:Introduction}

For simplicity, we demonstrate the idea by using the second order elliptic problem that seeks an
unknown function $u$ satisfying
\begin{eqnarray}
-\nabla\cdot (\nabla u)&=&f\quad \mbox{in}\;\Omega,\label{pde}\\
u&=&g\quad\mbox{on}\;\partial\Omega,\label{bc}
\end{eqnarray}
where $\Omega$ is a polytopal domain in $\mathbb{R}^2$.

The weak form of the problem (\ref{pde})-(\ref{bc}) is to find $u\in H^1(\Omega)$
such that $u=g$ on $\partial\Omega$ and satisfies
\begin{eqnarray}
( \nabla u,\nabla v)=(f,v)\quad \forall v\in
H_0^1(\Omega).\label{weakform}
\end{eqnarray}

The weak Galerkin finite element method is an effective and flexible numerical technique 
for
solving partial differential equations.
It is a natural extension of the standard Galerkin finite element method
where classical derivatives were substituted by weakly defined
derivatives on functions with discontinuity. The WG method was first introduced in
\cite{wy,wymix} and then has been applied to solve various PDEs such as  second order elliptic equations,
 biharmonic equations, Stokes equations, Navier-Stokes equations, Brinkman equations, parabolic equations, Helmholtz equation, convection dominant problems, hyperbolic equations, and Maxwell's equations
\cite{cz,hmy,lyzz,Lin2014,mwy-helm,mwy-brinkman,mwy-biha,mwy-soe,mwyz-maxwell, mwyz-interface,Shields,ww-div-curl,wy-stokes}.

The main idea of weak Galerkin finite
element methods is the use of weak functions and their corresponding
 weak derivatives. For the second order
elliptic equation, weak functions have the form of $v=\{v_0,v_b\}$
with $v=v_0$ inside of each element and $v=v_b$ on the boundary of
the element. Both $v_0$ and $v_b$ can be approximated by polynomials
in $P_\ell(T)$ and $P_s(e)$ respectively, where $T$ stands for an
element and $e$ the edge or face of $T$, $\ell$ and $s$ are
non-negative integers with possibly different values. Weak
derivatives are defined for weak functions in the sense of
distributions. Denote by $G_m(T)$ the vector space for weak gradient.
Typical choices for $G_m(T)$ are $[P_m(T)]^d$ or $RT_m(T)$. Various combination of $(P_\ell(T),P_s(e),G_m(T))$
leads to different  weak Galerkin methods tailored for
specific partial differential equations.

Weak Galerkin finite element methods have two forms  for the problem (\ref{pde})-(\ref{bc}).
One is its standard formulation \cite{mwy-soe,wy}:  find $u_h\in V_h$
such that $u_h=Q_bg$ on $\partial\Omega$ and satisfies
\begin{eqnarray}
(\nabla_w u_h,\nabla_w v)+s(u_h, v)=(f,v)\quad \forall v\in V_h^0,\label{wg0}
\end{eqnarray}
where $s(\cdot,\cdot)$ is a parameter independent stabilizer.
Another one is WG stabilizer free formulation \cite{aw, yz, yzz}: find $u_h\in V_h$
such that $u_h=Q_bg$ on $\partial\Omega$ and satisfies
\begin{eqnarray}
(\nabla_w u_h,\nabla_w v)=(f,v)\quad \forall v\in V_h^0.\label{wg1}
\end{eqnarray}
Removing  stabilizers  simplifies the formulations and
reduces programming complexity.
A stabilizer free WG method can be obtained by raising the degree of polynomial
$m$ for approximating weak gradient in the WG element $(P_\ell(T),P_s(e),G_m(T))$.

The purpose of this paper is to investigate the performance of different WG elements
computationally in the weak Galerkin finite element methods with or without stabilizers.
Like a periodic table, we provide 31 tables that are informative and clearly
demonstrate special properties of each WG element.
We don't have all the theoretical answers for many interesting phenomena shown in the tables and
we leave them for interesting readers.

While preparing this manuscript,  three papers are in the process to answer some questions from the
   numerical results in the tables.
We are close to prove theoretically that the WG element $(P_k(T), P_k(e),[P_{k+1}]^2)$  has two orders of supercloseness
   in both energy norm and $L^2$ norm on rectangular meshes,  shown in Table \ref{t001rec}.
It is proved in \cite{awyz} that the WG element $(P_k(T), P_{k+1}(e),[P_{k+1}]^2)$  has two orders of supercloseness
   in both energy norm and $L^2$ norm, on general triangular meshes in Table \ref{t011tr}.
Due to the bad behavior of the  WG element $(P_k(T), P_{k-1}(e),[P_{k+1}]^2)$  shown in Table \ref{t102rec} and \ref{t102tr},
   a new definition of the weak gradient is introduced in \cite{Ye-Zhang-k-1} so that the
   element can still converge in optimal order on general polytopal meshes.

The WG methods are designed for using discontinuous approximations  on general polytopal meshes. Due to  limited space, we only consider the finite element partitions including rectangles and triangles.

\section{Weak Galerkin Finite Element formulations}

Let ${\cal T}_h$ be a partition of the domain $\Omega$ consisting of
rectangles or triangles. Denote by ${\cal E}_h$
the set of all edges in ${\cal T}_h$, and let ${\cal
E}_h^0={\cal E}_h\backslash\partial\Omega$ be the set of all
interior edges or flat faces. For every element $T\in \T_h$, we
denote by $h_T$ its diameter and mesh size $h=\max_{T\in\T_h} h_T$
for ${\cal T}_h$. Let $P_k(T)$ consist all the polynomials defined on $T$ of degree less or equal to $k$.

\begin{thm1}
For $T\in\T_h$ and $\ell,s\ge 0$, define a local WG element $W_{\ell,s}(T)$ as,
\begin{equation}\label{l-space}
W_{\ell,s}(T)=\{v=\{v_0,v_b\}:\; \; v_0|_T\in P_{\ell}(T),\; v_b|_e\in P_{s}(e),\; e\subset\pT\}.
\end{equation}
\end{thm1}

\begin{thm1}
For any $v=\{v_0,v_b\}\in W_{\ell,s}(T)$, a weak gradient $\nabla_wv\in G_m(T)$ is defined as a unique solution of the following equation
\begin{equation}\label{d-d}
(\nabla_{w} v, \bq)_T = -(v_0,\nabla\cdot \bq)_T+ \langle v_b,
\bq\cdot\bn\rangle_{\partial T}\qquad \forall \bq\in G_m(T).
\end{equation}
\end{thm1}
A typical choice of $G_m(T)$ is $[P_m(T)]^d$, or $RT_m(T)$. Different
combinations of $(\ell,s,m)$ associated with a WG element $W_{\ell,s}(T)/G_m(T)$  leads to different weak Galerkin
finite element formulations. The weak gradient $\nabla_w$ defined in (\ref{d-d}) is  an approximation of $\nabla$ that is computed on each element $T$.

\smallskip

\begin{remark}
Please note that the space $G_m(T)$ is used to calculated weak gradient and  does not introduce additional degrees of freedom to the resulting linear system.
\end{remark}

\smallskip

\begin{thm1}
Define a WG finite element space $V_h$ associated with $\T_h$ as follows
\begin{equation}\label{vh}
V_h=\{v=\{v_0,v_b\}:\; \; v|_T\in W_{\ell,s}(T),\; \forall T\in\T_h  \}.
\end{equation}
\end{thm1}
We would like to emphasize that any function $v\in V_h$ has a single
value $v_b$ on each edge $e\in\E_h$. The subspace of $V_h$ consisting of functions with vanishing boundary value is denoted as $V_h^0$.

Let $Q_0$ and $Q_b$ be the two element-wise defined $L^2$ projections onto $P_\ell(T)$ and $P_s(e)$  on each $T\in\T_h$, respectively. Define $Q_hu=\{Q_0u,Q_bu\}\in V_h$. Let $\Q_h$ be the element-wise defined $L^2$ projection onto $G_m(T)$ on each element $T\in\T_h$.

For simplicity, we adopt the following notations,
\begin{eqnarray*}
(v,w)_{\T_h} &=&\sum_{T\in\T_h}(v,w)_T=\sum_{T\in\T_h}\int_T vw d\bx,\\
 \l v,w\r_{\partial\T_h}&=&\sum_{T\in\T_h} \l v,w\r_\pT=\sum_{T\in\T_h} \int_\pT vw ds.
\end{eqnarray*}

\begin{algorithm}
A numerical approximation for (\ref{pde})-(\ref{bc}) can be
obtained by seeking $u_h=\{u_0,u_b\}\in V_h$
satisfying $u_b=Q_bg$ on $\partial\Omega$ and  the following equation:
\begin{equation}\label{wg}
(\nabla_wu_h,\nabla_wv)_{\T_h}+s(u_h,v)=(f,\; v_0) \quad\forall v=\{v_0,v_b\}\in V_h^0,
\end{equation} where the stabilizer $s(\cdot,\cdot)$  is defined as
\begin{equation}\label{s1}
s(u_h,v)= \sum_{T\in\T_h}h_T^j\l Q_b u_0-u_b,Q_b v_0-v_b\r_\pT.
\end{equation}
Let $j=\infty$ in \eqref{s1},  we mean $s(u_h,v)=0$, i.e., we have the following stabilizer-free WG formulation,
\begin{equation}\label{wg2}
(\nabla_wu_h,\nabla_wv)_{\T_h}=(f,\; v_0) \quad\forall v=\{v_0,v_b\}\in V_h^0.
\end{equation}
\end{algorithm}

In the following sections, we will conduct extensive numerical tests to study the performance of different WG elements and record the results in 31 tables. In all the tables below, $j=\infty$ refers to the stabilizer free WG formulation (\ref{wg2}), where $j$ is defined in (\ref{s1}).

\section{The WG elements with $\ell= s$ on rectangular mesh}

\begin{figure}[H]
 \begin{center} \setlength\unitlength{1.25pt}
\begin{picture}(260,80)(0,0)
  \def\tr{\begin{picture}(20,20)(0,0)\put(0,0){\line(1,0){20}}\put(0,20){\line(1,0){20}}
          \put(0,0){\line(0,1){20}} \put(20,0){\line(0,1){20}}
      \end{picture}}
 {\setlength\unitlength{5pt}
 \multiput(0,0)(20,0){1}{\multiput(0,0)(0,20){1}{\tr}}}

  {\setlength\unitlength{2.5pt}
 \multiput(45,0)(20,0){2}{\multiput(0,0)(0,20){2}{\tr}}}

  \multiput(180,0)(20,0){4}{\multiput(0,0)(0,20){4}{\tr}}

 \end{picture}\end{center}
\caption{The first three level rectangular grids.}
\label{grid1}
\end{figure}
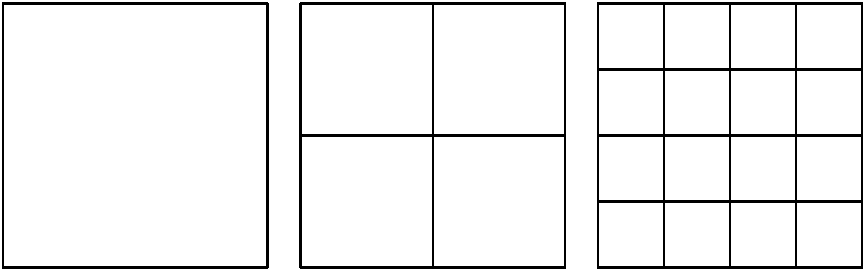

Next we will study convergence rate for the WG element $(P_\ell(T),P_s(e),G_m(T))$ with $\ell=s$ on rectangular meshes.
The rectangular meshes used in the computation are illustrated in
  Figure \ref{grid1}.

Table \ref{t110rec} demonstrates the convergence rates for  $(P_k(T),P_k(e),[P_{k-1}(T)]^2)$ with a stabilizer of different $j$ defined in (\ref{s1}) on rectangular mesh.

\begin{table}[H]
 \caption{Element $(P_k(T),P_k(e),[P_{k-1}(T)]^2)$ on rectangular mesh,
  $\3bar\cdot\3bar=O(h^{r_1})$ and $\|\cdot\|=O(h^{r_2})$.}
    \label{t110rec}
\begin{center}
\begin{tabular}{l|c|c|c|c|c|c|c}
\hline
element & $P_k(T)$ & $P_{k}(e)$ &  $[P_{k-1}(T)]^2$  & $j$ & $r_1$  & $r_2$ & Proved \\
\hline
\ref{t110rec}.1 & $ $ & $ $ & $ $ &  $-1$ &  $1 $ & $2$ & Yes \\
\ref{t110rec}.2 & $ $ & $ $ & $ $ &  $0$ & $0.5 $ & $1$ & No \\
\ref{t110rec}.3 & $P_1(T)$ & $P_1(e)$ & $[P_0(T)]^2$ &  $1$ & $0 $ & $0 $ & No \\
\ref{t110rec}.4 & $ $ & $ $ & $ $ & $\infty$ & $-\infty$ & $-\infty$& No   \\\hline
\ref{t110rec}.5 & $      $ & $      $ & $          $ & $-1$ &  $2 $ & $3$ & Yes \\
\ref{t110rec}.6 & $      $ & $      $ & $          $ & $0$ & $1.5$ & $2$ & No \\
\ref{t110rec}.7 & $P_2(T)$ & $P_2(e)$ & $[P_1(T)]^2$ &  $1$ & $1$ & $1$ & No \\
\ref{t110rec}.8 & $      $ & $      $ & $          $ & $\infty$ & $-\infty$ & $-\infty$& No   \\\hline
\ref{t110rec}.9 & $      $ & $      $  & $          $ & $-1$ &  $3 $ & $4$ & Yes \\
\ref{t110rec}.10 & $      $ & $      $ & $          $ & $0$ & $2.5$ & $3$ & No \\
\ref{t110rec}.11 & $P_3(T)$ & $P_3(e)$ & $[P_2(T)]^2$ & $1$ & $2$ & $2$ & No \\
\ref{t110rec}.12 & $      $ & $      $ & $          $ & $\infty$ & $-\infty$ & $-\infty$& No   \\\hline
\end{tabular}
\end{center}
\end{table}

Table \ref{t000rec}  demonstrates the convergence rates for  $(P_k(T),P_k(e),[P_{k}(T)]^2)$ with a stabilizer of different $j$ on rectangular mesh.

\begin{table}[H]
 \caption{Element $(P_k(T),P_k(e),[P_{k}(T)]^2)$ on rectangular mesh,
  $\3bar\cdot\3bar=O(h^{r_1})$ and $\|\cdot\|=O(h^{r_2})$.}
    \label{t000rec}
\begin{center}
\begin{tabular}{l|c|c|c|c|c|c|c}
\hline
element & $P_k(T)$ & $P_k(e)$ &  $[P_{k}(T)]^2$  & $j$ & $r_1$  & $r_2$ & Proved \\
\hline
\ref{t000rec}.1 & $      $ & $      $ & $          $ & $-1$ & $0 $ & $0 $ & No \\
\ref{t000rec}.2 & $      $ & $      $ & $          $ & $0$ & $0.5 $ & $1$ & No \\
\ref{t000rec}.3 & $P_0(T)$ & $P_0(e)$ & $[P_0(T)]^2$ & $1$ & $0 $ & $0 $ & No \\
\ref{t000rec}.4 & $      $ & $      $ & $          $ & $\infty$ & $-\infty$ & $-\infty$& No   \\\hline
\ref{t000rec}.5 & $      $ & $      $ & $          $ & $-1$  & $1$   & $2$ & No \\
\ref{t000rec}.6 & $      $ & $      $ & $          $ & $0$ & $1.5 $ & $2$ & No \\
\ref{t000rec}.7 & $P_1(T)$ & $P_1(e)$ & $[P_1(T)]^2$ & $1$ & $1$ & $1$ & No \\
\ref{t000rec}.8 & $      $ & $      $ & $          $ & $\infty$ & $-\infty$ & $-\infty$& No  \\\hline
\ref{t000rec}.9  & $      $ & $      $ & $          $ &  $-1$  & $2$   & $3$ & No \\
\ref{t000rec}.10 & $      $ & $      $ & $          $ & $0$ & $2.5 $ & $3$ & No \\
\ref{t000rec}.11 & $P_2(T)$ & $P_2(e)$ & $[P_2(T)]^2$ & $1$ & $2$ & $2$ & No \\
\ref{t000rec}.12 & $      $ & $      $ & $          $ & $\infty$ & $-\infty$ & $-\infty$& No   \\\hline
\end{tabular}
\end{center}
\end{table}

\begin{remark}
Theorem 4.9 in \cite{wwzz} guarantees the optimal convergence rate of the WG element \ref{t000rec}.5 in  the $\3bar\cdot\3bar$ norm. However the optimal convergence rate in the $L^2$ norm is not proved in \cite{wwzz}.
Therefore, we still mark $proved=No$ in the case \ref{t000rec}.5.
\end{remark}

\smallskip

Table \ref{t001rec}  demonstrates the convergence rates for  $(P_k(T),P_k(e),[P_{k+1}(T)]^2)$ with a stabilizer of different $j$ on rectangular mesh.

\begin{table}[H]
 \caption{Element $(P_k(T),P_k(e),[P_{k+1}(T)]^2)$ on rectangular mesh,
  $\3bar\cdot\3bar=O(h^{r_1})$ and $\|\cdot\|=O(h^{r_2})$.}
    \label{t001rec}
\begin{center}
\begin{tabular}{c|c|c|c|c|c|c|c}
\hline
element & $P_k(T)$ & $P_k(e)$ &  $[P_{k+1}(T)]^2$  & $j$ & $r_1$  & $r_2$ & Proved \\
\hline
\ref{t001rec}.1 & $      $ & $      $ & $          $ &  $-1$  & $0$   & $0$ & No \\
\ref{t001rec}.2 & $      $ & $      $ & $          $ &  $0$ & $1 $ & $1$ & No \\
\ref{t001rec}.3 & $P_0(T)$ & $P_0(e)$ & $[P_1(T)]^2$ &   $1$ & $2$ & $2$ & No \\
\ref{t001rec}.4 & $      $ & $      $ & $          $ & $\infty$ & $2$ & $2$& No   \\\hline
\hline
\ref{t001rec}.5 & $      $ & $      $ & $          $ &   $-1$  & $1$   & $2$ & No \\
\ref{t001rec}.6 & $      $ & $      $ & $          $ &  $0$ & $2 $ & $3$ & No \\
\ref{t001rec}.7 & $P_1(T)$ & $P_1(e)$ & $[P_2(T)]^2$ &  $1$ & $3$ & $4$ & No \\
\ref{t001rec}.8 & $      $ & $      $ & $          $ &  $\infty$ & $3$ & $4$& Yes   \\\hline
\hline
\ref{t001rec}.9  & $      $ & $      $ & $          $ &   $-1$  & $2$   & $3$ & No \\
\ref{t001rec}.10 & $      $ & $      $ & $          $ &  $0$ & $3$ & $4$ & No \\
\ref{t001rec}.11 & $P_2(T)$ & $P_2(e)$ & $[P_3(T)]^2$ &  $1$ & $4$ & $5$ & No \\
\ref{t001rec}.12 & $      $ & $      $ & $          $ & $\infty$ & $4$ & $5$& Yes   \\\hline
\hline
\end{tabular}
\end{center}
\end{table}

Table \ref{t002rec}   demonstrates the convergence rates for  $(P_k(T),P_k(e),[P_{k+2}(T)]^2)$ with a stabilizer of different $j$ on rectangular mesh.

\begin{table}[H]
 \caption{Element $(P_k(T),P_k(e),[P_{k+2}(T)]^2)$ on rectangular mesh,
$\3bar\cdot\3bar=O(h^{r_1})$ and $\|\cdot\|=O(h^{r_2})$.}
    \label{t002rec}
\begin{center}
\begin{tabular}{c|c|c|c|c|c|c|c}
\hline
element & $P_k(T)$ & $P_k(e)$ &  $[P_{k+2}(T)]^2$  & $j$ & $r_1$  & $r_2$ & Proved \\
\hline
\ref{t002rec}.1 & $      $ & $      $ & $          $ &  $-1$  & $0$   & $0$ & No \\
\ref{t002rec}.2 & $      $ & $      $ & $          $ &  $0$ & $0 $ & $0$ & No \\
\ref{t002rec}.3 & $P_0(T)$ & $P_0(e)$ & $[P_2(T)]^2$ &   $1$ & $0$ & $0$ & No \\
\ref{t002rec}.4 & $      $ & $      $ & $          $ & $\infty$ & $0$ & $0$& No   \\\hline
\hline
\ref{t002rec}.5 & $      $ & $      $ & $          $ &   $-1$  & $1$   & $2$ & No \\
\ref{t002rec}.6 & $      $ & $      $ & $          $ &  $0$ & $1 $ & $2$ & No \\
\ref{t002rec}.7 & $P_1(T)$ & $P_1(e)$ & $[P_3(T)]^2$ &  $1$ & $1$ & $2$ & No \\
\ref{t002rec}.8 & $      $ & $      $ & $          $ &  $\infty$ & $1$ & $2$& Yes   \\\hline
\hline
\ref{t002rec}.9  & $      $ & $      $ & $          $ &   $-1$  & $2$   & $3$ & No \\
\ref{t002rec}.10 & $      $ & $      $ & $          $ &  $0$ & $ 2 $ & $3$ & No \\
\ref{t002rec}.11 & $P_2(T)$ & $P_2(e)$ & $[P_4(T)]^2$ &  $1$ &   $2$ & $3$ & No \\
\ref{t002rec}.12 & $      $ & $      $ & $          $ & $\infty$ & $2$ & $3$& Yes  \\\hline
\hline
\end{tabular}
\end{center}
\end{table}

\begin{remark}
For the $P_k(T)-P_k(e)$ element, Tables \ref{t110rec}-\ref{t002rec} demonstrate that the performance of the WG solutions are getting better when the degree of the polynomials for weak gradient is increasing from $k-1$ to $k+1$. Specially the WG element $(P_k(T), P_k(e), [P_{k+1}(T)]^2)$ shows order two supercloseness in Table \ref{t001rec}. However, for the element $(P_k(T), P_k(e), [P_{k+2}(T)]^2$, the numerical tests in Table \ref{t002rec} show the convergence rate of the WG solution decreasing. Remember that increasing $m$ in $[P_m(T)]^2$ for weak gradient does not introduce additional degrees of freedom for the resulting linear systems.
\end{remark}

\smallskip

Table \ref{t00RT0rec}   demonstrates the convergence rates for  $(P_k(T),P_k(e), RT_k(T))$ with a stabilizer of different $j$ on rectangular mesh.

\begin{table}[H]
 \caption{Element $(P_k(T),P_k(e),RT_k(T))$ on rectangular mesh,
$\3bar\cdot\3bar=O(h^{r_1})$ and $\|\cdot\|=O(h^{r_2})$.}
    \label{t00RT0rec}
\begin{center}
\begin{tabular}{c|c|c|c|c|c|c|c}
\hline
element & $P_k(T)$ & $P_k(e)$ &  $RT_k(T)$  & $j$ & $r_1$  & $r_2$ & Proved \\
\hline
\ref{t00RT0rec}.1 & $      $ & $      $ & $          $ & $-1$ & $0$ & $0$ & No \\
\ref{t00RT0rec}.2 & $      $ & $      $ & $          $ & 0  & $1$ & $1$ & No \\
\ref{t00RT0rec}.3 & $P_0(T)$ & $P_0(e)$ & $RT_0(T)$ & 1 & $ 2$ & $2$ & No \\
\ref{t00RT0rec}.4 & $      $ & $      $ & $          $ &$\infty$& $2$ & $2$ & No \\
\hline
\ref{t00RT0rec}.5& $      $ & $      $ & $          $ & $ -1$ & $1$ & $2 $ & No  \\
\ref{t00RT0rec}.6& $      $ & $      $ & $          $ & 0 & $1.5$ & $2$ & No  \\
\ref{t00RT0rec}.7 & $P_1(T)$ & $P_1(e)$ & $RT_1(T)$ & 1 & $1$ & $1$ & No  \\
\ref{t00RT0rec}.8& $      $ & $      $ & $          $ &$\infty$& $1$ & $1$ & No  \\
\hline
\ref{t00RT0rec}.9 & $      $ & $      $ & $          $ & $-1$ & $2$ & $3 $ & No  \\
\ref{t00RT0rec}.10& $      $ & $      $ & $          $ & 0  & $2.5$ & $3$ & No  \\
\ref{t00RT0rec}.11 & $P_2(T)$ & $P_2(e)$ & $RT_2(T)$ & 1 & $2$ & $2$ & No  \\
\ref{t00RT0rec}.12& $      $ & $      $ & $          $ &$\infty$& $2$ & $2$ & No \\
\hline
\end{tabular}
\end{center}
\end{table}

\section{The WG elements for $\ell= s$ on triangular mesh}

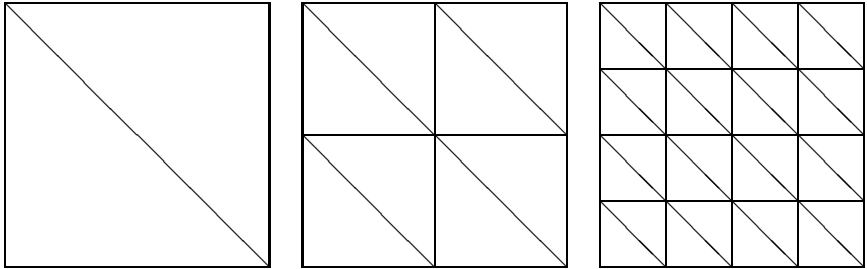
\begin{figure}[h!]
 \begin{center} \setlength\unitlength{1.25pt}
\begin{picture}(260,80)(0,0)
  \def\tr{\begin{picture}(20,20)(0,0)\put(0,0){\line(1,0){20}}\put(0,20){\line(1,0){20}}
          \put(0,0){\line(0,1){20}} \put(20,0){\line(0,1){20}}
   \put(0,20){\line(1,-1){20}}   \end{picture}}
 {\setlength\unitlength{5pt}
 \multiput(0,0)(20,0){1}{\multiput(0,0)(0,20){1}{\tr}}}

  {\setlength\unitlength{2.5pt}
 \multiput(45,0)(20,0){2}{\multiput(0,0)(0,20){2}{\tr}}}

  \multiput(180,0)(20,0){4}{\multiput(0,0)(0,20){4}{\tr}}

 \end{picture}\end{center}
\caption{The first three level triangular meshes.}
\label{grid2}
\end{figure}

The triangular meshes used in the computation are displayed in
Figure \ref{grid2}.

Table \ref{t110tr} demonstrates the convergence rates for  $(P_k(T),P_k(e),[P_{k-1}(T)]^2)$ with a stabilizer of different $j$ on triangular mesh.

\begin{table}[H]
 \caption{Element $(P_k(T),P_k(e),[P_{k-1}(T)]^2)$ on triangular mesh,
  $\3bar\cdot\3bar=O(h^{r_1})$ and $\|\cdot\|=O(h^{r_2})$.}
    \label{t110tr}
\begin{center}
\begin{tabular}{c|c|c|c|c|c|c|c}
\hline
element & $P_k(T)$ & $P_k(e)$ &  $[P_{k-1}(T)]^2$  & $j$ & $r_1$  & $r_2$ & Proved \\
\hline
\ref{t110tr}.1 & $      $ & $      $ & $          $ & $-1$ & $1$ & $2$ &Yes \\
\ref{t110tr}.2 & $      $ & $      $ & $          $ & 0  & $0.5$ & $1$  & No \\
\ref{t110tr}.3 & $P_1(T)$ & $P_1(e)$ & $[P_0(T)]^2$ & 1 & $0$ & $0$ & No  \\
\ref{t110tr}.4 & $      $ & $      $ & $          $ &$\infty$& $-\infty$ & $-\infty$&No  \\
\hline
\ref{t110tr}.5 & $      $ & $      $ & $          $ & $-1$ & $2$ & $3$ &Yes   \\
\ref{t110tr}.6 & $      $ & $      $ & $          $ & 0 & $1.5$ & $2$&No   \\
\ref{t110tr}.7 & $P_2(T)$ & $P_2(e)$ & $[P_1(T)]^2$ & 1 & $1$ & $1$  &No  \\
\ref{t110tr}.8 & $      $ & $      $ & $          $ &$\infty$& $-\infty$ & $-\infty$&No  \\
\hline
\ref{t110tr}.9  & $      $ & $      $ & $          $ & $-1$ & $3$ & $4$ &Yes     \\
\ref{t110tr}.10 & $      $ & $      $ & $          $ & 0 & $2.5$ & $3$  &No  \\
\ref{t110tr}.11 & $P_3(T)$ & $P_3(e)$ & $[P_2(T)]^2$ & 1 & $2$ & $2$  &No  \\
\ref{t110tr}.12 & $      $ & $      $ & $          $ &$\infty$& $-\infty$ & $-\infty$&No  \\
\hline
\end{tabular}
\end{center}
\end{table}

Table \ref{t000tr} demonstrates the convergence rates for  $(P_k(T),P_k(e),[P_{k}(T)]^2)$ with a stabilizer of different $j$ on triangular mesh.

\begin{table}[H]
 \caption{Element $(P_k(T),P_k(e),[P_{k}(T)]^2)$ on triangular mesh,
  $\3bar\cdot\3bar=O(h^{r_1})$ and $\|\cdot\|=O(h^{r_2})$.}
    \label{t000tr}
\begin{center}
\begin{tabular}{c|c|c|c|c|c|c|c}
\hline
element & $P_k(T)$ & $P_k(e)$ &  $[P_{k}(T)]^2$  & $j$ & $r_1$  & $r_2$ & Proved \\
\hline
\ref{t000tr}.1 & $      $ & $      $ & $          $ & $-1$ & $0$ & $0$ &No \\
\ref{t000tr}.2 & $      $ & $      $ & $          $ & 0  & $0.5$ & $1$  & No \\
\ref{t000tr}.3 & $P_0(T)$ & $P_0(e)$ & $[P_0(T)]^2$ & 1 & $0$ & $0$ & No  \\
\ref{t000tr}.4 & $      $ & $      $ & $          $ &$\infty$& $-\infty$ & $-\infty$&No  \\
\hline
\ref{t000tr}.5 & $      $ & $      $ & $          $ & $-1$ & $1$ & $2$ &No   \\
\ref{t000tr}.6 & $      $ & $      $ & $          $ & 0 & $1.5$ & $2$&No   \\
\ref{t000tr}.7 & $P_1(T)$ & $P_1(e)$ & $[P_1(T)]^2$ & 1 & $1$ & $1$  &No  \\
\ref{t000tr}.8 & $      $ & $      $ & $          $ &$\infty$& $-\infty$ & $-\infty$&No  \\
\hline
\ref{t000tr}.9  & $      $ & $      $ & $          $ & $-1$ & $2$ & $3$ &No  \\
\ref{t000tr}.10 & $      $ & $      $ & $          $ & 0 & $2.5$ & $3$  &No  \\
\ref{t000tr}.11 & $P_2(T)$ & $P_2(e)$ & $[P_2(T)]^2$ & 1 & $2$ & $2$  &No  \\
\ref{t000tr}.12 & $      $ & $      $ & $          $ &$\infty$& $-\infty$ & $-\infty$&No  \\
\hline
\end{tabular}
\end{center}
\end{table}

Table \ref{t001tr}  demonstrates the convergence rates for  $(P_k(T),P_k(e),[P_{k+1}(T)]^2)$ with a stabilizer of different $j$ on triangular mesh.

\begin{table}[H]
 \caption{Element $(P_k(T),P_k(e),[P_{k+1}(T)]^2)$ on triangular mesh,
  $\3bar\cdot\3bar=O(h^{r_1})$ and $\|\cdot\|=O(h^{r_2})$.}
    \label{t001tr}
\begin{center}
\begin{tabular}{c|c|c|c|c|c|c|c}
\hline
element & $P_k(T)$ & $P_k(e)$ &  $[P_{k+1}(T)]^2$  & $j$ & $r_1$  & $r_2$ & Proved \\
\hline
\ref{t001tr}.1 & $      $ & $      $ & $          $ & $-1$ & $0$ & $0$ &No \\
\ref{t001tr}.2 & $      $ & $      $ & $          $ & 0  & $0 $ & $0$  & No \\
\ref{t001tr}.3 & $P_0(T)$ & $P_0(e)$ & $[P_1(T)]^2$ & 1 & $0$ & $0$ & No  \\
\ref{t001tr}.4 & $      $ & $      $ & $          $ &$\infty$& $0$ & $0$&No  \\
\hline
\ref{t001tr}.5 & $      $ & $      $ & $          $ & $-1$ & $1$ & $2$ &No   \\
\ref{t001tr}.6 & $      $ & $      $ & $          $ & 0 & $1$ & $2$&No   \\
\ref{t001tr}.7 & $P_1(T)$ & $P_1(e)$ & $[P_2(T)]^2$ & 1 & $1$ & $2$  &No  \\
\ref{t001tr}.8 & $      $ & $      $ & $          $ &$\infty$& $1$ & $2$&Yes  \\
\hline
\ref{t001tr}.9  & $      $ & $      $ & $          $ & $-1$ & $2$ & $3$ &No  \\
\ref{t001tr}.10 & $      $ & $      $ & $          $ & 0 & $2 $ & $3$  &No  \\
\ref{t001tr}.11 & $P_2(T)$ & $P_2(e)$ & $[P_3(T)]^2$ & 1 & $2$ & $3$  &No  \\
\ref{t001tr}.12 & $      $ & $      $ & $          $ &$\infty$& $2$ & $3$& Yes  \\
\hline
\end{tabular}
\end{center}
\end{table}

\begin{remark}
The WG element $(P_k(T), P_k(e), [P_{k+1}(T)]^2)$ performs much better on rectangular meshes than triangular meshes.
\end{remark}

\smallskip

Table \ref{t00RTtr}  demonstrates the convergence rates for  $(P_k(T),P_k(e), RT_k(T))$ with a stabilizer of different $j$ on triangular mesh.

\begin{table}[H]
 \caption{Element $(P_k(T),P_k(e),RT_k(T))$ on triangular mesh,
  $\3bar\cdot\3bar=O(h^{r_1})$ and $\|\cdot\|=O(h^{r_2})$.}
    \label{t00RTtr}
\begin{center}
\begin{tabular}{c|c|c|c|c|c|c|c}
\hline
element & $P_k(T)$ & $P_k(e)$ &  $RT_k(T)$  & $j$ & $r_1$  & $r_2$ & Proved \\
\hline
\ref{t00RTtr}.1 & $      $ & $      $ & $          $ & $-1$ & $0$ & $0$ &No \\
\ref{t00RTtr}.2 & $      $ & $      $ & $          $ & 0  & $1 $ & $1$  & No \\
\ref{t00RTtr}.3 & $P_0(T)$ & $P_0(e)$ & $RT_0(T)$ & 1 & $1$ & $2$ & No  \\
\ref{t00RTtr}.4 & $      $ & $      $ & $          $ &$\infty$& $1$ & $2$&Yes  \\
\hline
\ref{t00RTtr}.5 & $      $ & $      $ & $          $ & $-1$ & $1$ & $2$ & No  \\
\ref{t00RTtr}.6 & $      $ & $      $ & $          $ & 0 & $2$ & $3$&No   \\
\ref{t00RTtr}.7 & $P_1(T)$ & $P_1(e)$ & $RT_1(T)$ & 1 & $2$ & $3$  &No  \\
\ref{t00RTtr}.8 & $      $ & $      $ & $          $ &$\infty$& $2$ & $3$& Yes  \\
\hline
\ref{t00RTtr}.9  & $      $ & $      $ & $          $ & $-1$ & $2$ & $3$ &No \\
\ref{t00RTtr}.10 & $      $ & $      $ & $          $ & 0 & $3 $ & $4$  &No  \\
\ref{t00RTtr}.11 & $P_2(T)$ & $P_2(e)$ & $RT_2(T)$ & 1 & $3$ & $4$  &No  \\
\ref{t00RTtr}.12 & $      $ & $      $ & $          $ &$\infty$& $3$ & $4$& Yes \\
\hline
\end{tabular}
\end{center}
\end{table}

\section{The WG elements with $\ell< s$ on  rectangular mesh}

The following  table demonstrates the convergence rates for  $(P_k(T),P_{k+1}(e),[P_{k-1}(T)]^2)$ with a stabilizer of different $j$ on rectangular mesh.

\begin{table}[H]
 \caption{Element $(P_k(T),P_{k+1}(e),[P_{k-1}(T)]^2)$ on rectangular mesh,
  $\3bar\cdot\3bar=O(h^{r_1})$ and $\|\cdot\|=O(h^{r_2})$.}
    \label{t120rec}
\begin{center}
\begin{tabular}{c|c|c|c|c|c|c|c}
\hline
element & $P_k(T)$ & $P_{k+1}(e)$ &  $[P_{k-1}(T)]^d$  & $j$ & $r_1$  & $r_2$ & Proved \\
\hline
\ref{t120rec}.1 & $      $ & $      $ & $          $ &  $-1$ & $1$ & $2$ & No   \\
\ref{t120rec}.2 & $      $ & $      $ & $          $ &  $0 $ & $0.5$ & $1$ & No   \\
\ref{t120rec}.3 & $P_1(T)$ & $P_2(e)$ & $[P_0(T)]^2$ & $1 $ & $ 0$ & $0$ & No   \\
\ref{t120rec}.4 & $      $ & $      $ & $          $ &$\infty$ & $-\infty$ & $-\infty$& No   \\
\hline
\ref{t120rec}.5 & $      $ & $      $ & $          $ & $-1$ & $2$ & $3$  &  No   \\
\ref{t120rec}.6 & $      $ & $      $ & $          $ & 0 & $1.5$ & $2$  & No   \\
\ref{t120rec}.7 & $P_2(T)$ & $P_3(e)$ & $[P_1(T)]^2$ &$1 $ & $1$ & $1$  & No   \\
\ref{t120rec}.8 & $      $ & $      $ & $          $ &$\infty$ & $-\infty$ & $-\infty$& No   \\
\hline
\ref{t120rec}.9 & $      $ & $      $ & $          $ & $-1$ &  $3 $ & $4$  &  No   \\
\ref{t120rec}.10 & $      $ & $      $ & $          $ & 0 &  $2.5$ & $3$  & No   \\
\ref{t120rec}.11 & $P_3(T)$ & $P_4(e)$ & $[P_2(T)]^2$ &$1 $ &$2$ & $2$  & No   \\
\ref{t120rec}.12 & $      $ & $      $ & $          $ &$\infty$ & $-\infty$ & $-\infty$& No   \\
\hline
\end{tabular}
\end{center}
\end{table}

The following  table demonstrates the convergence rates for  $(P_k(T),P_{k+1}(e),[P_{k}(T)]^2)$ with a stabilizer of different $j$ on rectangular mesh.

\begin{table}[H]
 \caption{Element $(P_k(T),P_{k+1}(e),[P_{k}(T)]^2)$ on rectangular mesh,
  $\3bar\cdot\3bar=O(h^{r_1})$ and $\|\cdot\|=O(h^{r_2})$.}
    \label{t010rec}
\begin{center}
\begin{tabular}{c|c|c|c|c|c|c|c}
\hline
element & $P_k(T)$ & $P_{k+1}(e)$ &  $[P_{k}(T)]^2$  & $j$ & $r_1$  & $r_2$ & Proved \\
\hline
\ref{t010rec}.1 & $      $ & $      $ & $          $ & $-1$ & $0$ & $0$ & No   \\
\ref{t010rec}.2 & $      $ & $      $ & $          $ & $0 $ & $0.5$ & $1$ & No   \\
\ref{t010rec}.3 & $P_0(T)$ & $P_1(e)$ & $[P_0(T)]^2$ & $1 $ & $0$ & $0$ & No   \\
\ref{t010rec}.4 & $      $ & $      $ & $          $ & $\infty$ & $-\infty$ & $-\infty$& No   \\
\hline
\ref{t010rec}.5 & $      $ & $      $ & $          $ & $-1$ & $1$ & $2$ &  No   \\
\ref{t010rec}.6 & $      $ & $      $ & $          $ & $0 $ & $1.5$ & $2$  & No   \\
\ref{t010rec}.7 & $P_1(T)$ & $P_2(e)$ & $[P_1(T)]^2$ & $1 $ &$1$ & $1$  & No   \\
\ref{t010rec}.8 & $      $ & $      $ & $          $ & $\infty$ & $-\infty$ & $-\infty$& No   \\
\hline
\ref{t010rec}.9  & $      $ & $      $ & $          $ & $-1$ & $2$ & $3$ &  No   \\
\ref{t010rec}.10 & $      $ & $      $ & $          $ & $0 $ & $2.5$ & $3$  & No   \\
\ref{t010rec}.11 & $P_2(T)$ & $P_3(e)$ & $[P_2(T)]^2$ & $1 $ & $2$ & $2$  & No   \\
\ref{t010rec}.12 & $      $ & $      $ & $          $ & $\infty$ & $-\infty$ & $-\infty$& No   \\
\hline
\end{tabular}
\end{center}
\end{table}

Table \ref{t011rec} demonstrates the convergence rates for  $(P_k(T),P_{k+1}(e),[P_{k+1}(T)]^2)$ with a stabilizer of different $j$ on rectangular mesh.

\begin{table}[H]
 \caption{Element $(P_k(T),P_{k+1}(e),[P_{k+1}(T)]^2)$ on rectangular mesh,
  $\3bar\cdot\3bar=O(h^{r_1})$ and $\|\cdot\|=O(h^{r_2})$.}
    \label{t011rec}
\begin{center}
\begin{tabular}{c|c|c|c|c|c|c|c}
\hline
element & $P_k(T)$ & $P_{k+1}(e)$ &  $[P_{k+1}(T)]^2$  & $j$ & $r_1$  & $r_2$ & Proved \\
\hline
\ref{t011rec}.1 & $      $ & $      $ & $          $ & $-1$ & $0$ & $0$& No   \\
\ref{t011rec}.2 & $      $ & $      $ & $          $ & $0 $ & $0.5$ & $1$  & No \\
\ref{t011rec}.3 & $P_0(T)$ & $P_1(e)$ & $[P_1(T)]^2$ & $1 $& $1$ & $2$ & No  \\
\ref{t011rec}.4 & $      $ & $      $ & $          $ & $\infty$& $1$ & $2$ & No  \\
\hline
\ref{t011rec}.5 & $      $ & $      $ & $          $ & $-1$ & $1$ & $2$ &  No  \\
\ref{t011rec}.6 & $      $ & $      $ & $          $ & $0 $& $ 1.5$ & $3$ & No  \\
\ref{t011rec}.7 & $P_1(T)$ & $P_2(e)$ & $[P_2(T)]^2$ & $1 $& $ 2 $ & $4$ & No     \\
\ref{t011rec}.8 & $      $ & $      $ & $          $ &$\infty$& $ 2 $ & $4$ & No     \\
\hline
\ref{t011rec}.9  & $      $ & $      $ & $          $ &$ -1 $& $2$ & $3$  &  No     \\
\ref{t011rec}.10 & $      $ & $      $ & $          $ &$ 0 $& $2.5$ & $4$   & No    \\
\ref{t011rec}.11 & $P_2(T)$ & $P_3(e)$ & $[P_3(T)]^2$ &$ 1 $& $3$ & $5$   & No    \\
\ref{t011rec}.12 & $      $ & $      $ & $          $ &$\infty$& $3$ & $5$  & No     \\
\hline
\end{tabular}
\end{center}
\end{table}

Table \ref{t012rec}  demonstrates the convergence rates for  $(P_k(T),P_{k+1}(e),[P_{k+2}(T)]^2)$ with a stabilizer of different $j$ on rectangular mesh.

\begin{table}[H]
 \caption{Element $(P_k(T),P_{k+1}(e),[P_{k+2}(T)]^2)$ on rectangular mesh,
  $\3bar\cdot\3bar=O(h^{r_1})$ and $\|\cdot\|=O(h^{r_2})$.}
    \label{t012rec}
\begin{center}
\begin{tabular}{c|c|c|c|c|c|c|c}
\hline
element & $P_k(T)$ & $P_{k+1}(e)$ &  $[P_{k+2}(T)]^2$  & $j$ & $r_1$  & $r_2$ & Proved \\
\hline
\ref{t012rec}.1 & $      $ & $      $ & $          $ &$ -1 $ & $0$ & $0$  & No  \\
\ref{t012rec}.2 & $      $ & $      $ & $          $ &$ 0  $ & $0$ & $0$  & No  \\
\ref{t012rec}.3 & $P_0(T)$ & $P_1(e)$ & $[P_2(T)]^2$ &$ 1 $   & $0$ & $0$  & No  \\
\ref{t012rec}.4 & $      $ & $      $ & $          $ &$\infty$ & $0$ & $0$  & No  \\
\hline
\ref{t012rec}.5 & $      $ & $      $ & $          $ &$ -1 $& $1$ & $2$  &  No  \\
\ref{t012rec}.6 & $      $ & $      $ & $          $ &$ 0 $& $1$ & $2$   & No   \\
\ref{t012rec}.7 & $P_1(T)$ & $P_2(e)$ & $[P_3(T)]^2$ &$ 1 $& $  1$ & $2  $   & No  \\
\ref{t012rec}.8 & $      $ & $      $ & $          $ &$\infty$& $1$ & $2$  & No \\
\hline
\ref{t012rec}.9  & $      $ & $      $ & $          $ &$ -1 $& $2$ & $3$  &  No  \\
\ref{t012rec}.10 & $      $ & $      $ & $          $ &$ 0 $& $2$ & $3$   & No   \\
\ref{t012rec}.11 & $P_2(T)$ & $P_3(e)$ & $[P_4(T)]^2$ &$ 1 $& $2$ & $3$   & No   \\
\ref{t012rec}.12 & $      $ & $      $ & $          $ &$\infty$& $2$ & $3$   & No  \\
\hline
\end{tabular}
\end{center}
\end{table}

\begin{remark}
The WG element $(P_k(T), P_{k+1}(e), [P_{k+1}(T)]^2)$ performs  better than WG element $(P_k(T), P_{k+1}(e), [P_{k+2}(T)]^2)$  although the later element uses higher degree polynomials for weak gradient.
\end{remark}

\smallskip

Table \ref{t01RTrec}  demonstrates the convergence rates for
 $(P_k(T),P_{k+1}(e), RT_k(T))$ with a stabilizer of different $j$ on rectangular mesh.

\begin{table}[H]
 \caption{Element $(P_k(T),P_{k+1}(e),RT_k(T))$ on rectangular mesh,
  $\3bar\cdot\3bar=O(h^{r_1})$ and $\|\cdot\|=O(h^{r_2})$.}
    \label{t01RTrec}
\begin{center}
\begin{tabular}{c|c|c|c|c|c|c|c}
\hline
element & $P_k(T)$ & $P_{k+1}(e)$ &  $RT_k(T)$  & $j$ & $r_1$  & $r_2$ & Proved \\
\hline
\ref{t01RTrec}.1& $      $ & $      $ & $          $ &$ -1 $ & $0$ & $0$  & No  \\
\ref{t01RTrec}.2& $      $ & $      $ & $          $ &$ 0  $ & $0.5$ & $1$  & No  \\
\ref{t01RTrec}.3 & $P_0(T)$ & $P_1(e)$ & $RT_0(T)$ &$ 1 $   & $1$ & $2$  & No  \\
\ref{t01RTrec}.4& $      $ & $      $ & $          $ &$\infty$ & $2$ & $2$  & No  \\
\hline
\ref{t01RTrec}.5& $      $ & $      $ & $          $ &$ -1 $& $1$ & $2$  & No  \\
\ref{t01RTrec}.6& $      $ & $      $ & $          $ &$ 0 $& $1.5$ & $2$   & No   \\
\ref{t01RTrec}.7 & $P_1(T)$ & $P_2(e)$ & $RT_1(T)$ &$ 1 $& $  1$ & $1$   & No  \\
\ref{t01RTrec}.8& $      $ & $      $ & $          $ &$\infty$& $-\infty$ & $-\infty$  &No \\
\hline
\ref{t01RTrec}.9 & $      $ & $      $ & $          $ &$ -1 $& $2$ & $3$  & No  \\
\ref{t01RTrec}.10& $      $ & $      $ & $          $ &$ 0 $& $2.5$ & $3$   & No   \\
\ref{t01RTrec}.11& $P_2(T)$ & $P_3(e)$ & $RT_2(T)$ &$ 1 $& $2$ & $2$   & No   \\
\ref{t01RTrec}.12& $      $ & $      $ & $          $ &$\infty$& $-\infty$ & $-\infty$  &No \\
\hline
\end{tabular}
\end{center}
\end{table}

Table \ref{t01RT1rec}  demonstrates the convergence rates for
 $(P_k(T),P_{k+1}(e), RT_{k+1}(T))$ with a stabilizer of different $j$ on rectangular mesh.

\begin{table}[H]
 \caption{Element $(P_k(T),P_{k+1}(e),RT_{k+1}(T))$ on rectangular mesh,
  $\3bar\cdot\3bar=O(h^{r_1})$ and $\|\cdot\|=O(h^{r_2})$.}
    \label{t01RT1rec}
\begin{center}
\begin{tabular}{c|c|c|c|c|c|c|c}
\hline
element & $P_k(T)$ & $P_{k+1}(e)$ &  $RT_{k+1}(T)$  & $j$ & $r_1$  & $r_2$ & Proved \\
\hline
\ref{t01RT1rec}.1 & $      $ & $      $ & $          $ &$ -1 $ & $0$ & $0$  & No  \\
\ref{t01RT1rec}.2 & $      $ & $      $ & $          $ &$ 0  $ & $0$ & $1$  & No  \\
\ref{t01RT1rec}.3 & $P_0(T)$ & $P_1(e)$ & $RT_1(T)$ &$ 1 $   & $0$ & $2$  & No  \\
\ref{t01RT1rec}.4 & $      $ & $      $ & $          $ &$\infty$ & $0$ & $2$  & No  \\
\hline
\ref{t01RT1rec}.5 & $      $ & $      $ & $          $ &$ -1 $& $1$ & $2$  & No  \\
\ref{t01RT1rec}.6 & $      $ & $      $ & $          $ &$ 0 $& $ 1$ & $2$   & No   \\
\ref{t01RT1rec}.7 & $P_1(T)$ & $P_2(e)$ & $RT_2(T)$ &$ 1 $& $  1$ & $2$   & No  \\
\ref{t01RT1rec}.8 & $      $ & $      $ & $          $ &$\infty$& $1$ & $2$  &No \\
\hline
\ref{t01RT1rec}.9 & $      $ & $      $ & $          $ &$ -1 $& $2$ & $4$  & No  \\
\ref{t01RT1rec}.10& $      $ & $      $ & $          $ &$ 0 $& $2$ & $4$   & No   \\
\ref{t01RT1rec}.11& $P_2(T)$ & $P_3(e)$ & $RT_3(T)$ &$ 1 $& $2$ & $4$   & No   \\
\ref{t01RT1rec}.12& $      $ & $      $ & $          $ &$\infty$& $2$ & $4$  &No \\
\hline
\end{tabular}
\end{center}
\end{table}

\section{The WG elements for $\ell< s$ on triangular mesh}

Table \ref{t120tr}   demonstrates the convergence rates for  $(P_k(T),P_{k+1}(e),[P_{k-1}(T)]^2)$ with a stabilizer of different $j$ on triangular mesh.

\begin{table}[H]
 \caption{Element $(P_k(T),P_{k+1}(e),[P_{k-1}(T)]^2)$ on triangular mesh,
  $\3bar\cdot\3bar=O(h^{r_1})$ and $\|\cdot\|=O(h^{r_2})$.}
    \label{t120tr}
\begin{center}
\begin{tabular}{c|c|c|c|c|c|c|c}
\hline
element & $P_k(T)$ & $P_{k+1}(e)$ &  $[P_{k-1}(T)]^2$  & $j$ & $r_1$  & $r_2$ & Proved \\
\hline
\ref{t120tr}.1 & $      $ & $      $ & $          $ &$ -1 $ & $1$ & $2$  &     No  \\
\ref{t120tr}.2 & $      $ & $      $ & $          $ &$ 0  $ & $0.5 $ & $1$  & No  \\
\ref{t120tr}.3 & $P_1(T)$ & $P_2(e)$ & $[P_{0}(T)]^2$ &$ 1 $   & $0$ & $0$  & No  \\
\ref{t120tr}.4 & $      $ & $      $ & $          $ &$\infty$ &$-\infty$ & $-\infty$  &No \\
\hline
\ref{t120tr}.5 & $      $ & $      $ & $          $ &$ -1 $& $2$ & $3$  &  No  \\
\ref{t120tr}.6 & $      $ & $      $ & $          $ &$ 0 $& $1.5$ & $2$   & No   \\
\ref{t120tr}.7 & $P_2(T)$ & $P_3(e)$ & $[P_{1}(T)]^2$ &$ 1 $& $  1$ & $1$   & No  \\
\ref{t120tr}.8 & $      $ & $      $ & $          $ &$\infty$& $-\infty$ & $-\infty$  &No \\
\hline
\ref{t120tr}.9 & $      $ & $      $ & $          $ &$ -1 $& $3$ & $4$  &  No  \\
\ref{t120tr}.10& $      $ & $      $ & $          $ &$ 0 $& $2.5$ & $3$   & No   \\
\ref{t120tr}.11& $P_3(T)$ & $P_4(e)$ & $[P_2(T)]^2$ &$ 1 $& $2$ & $2$   & No   \\
\ref{t120tr}.12& $      $ & $      $ & $          $ &$\infty$& $-\infty$ & $-\infty$  &No \\
\hline
\end{tabular}
\end{center}
\end{table}

Table \ref{t010tr}  demonstrates the convergence rates for  $(P_k(T),P_{k+1}(e),[P_{k}(T)]^2)$ with a stabilizer of different $j$ on triangular mesh.

\begin{table}[H]
 \caption{Element $(P_k(T),P_{k+1}(e),[P_{k}(T)]^2)$ on triangular mesh,
  $\3bar\cdot\3bar=O(h^{r_1})$ and $\|\cdot\|=O(h^{r_2})$.}
    \label{t010tr}
\begin{center}
\begin{tabular}{c|c|c|c|c|c|c|c}
\hline
element & $P_k(T)$ & $P_{k+1}(e)$ &  $[P_{k}(T)]^2$  & $j$ & $r_1$  & $r_2$ & Proved \\
\hline
\ref{t010tr}.1 & $      $ & $      $ & $          $ &$ -1 $ & $0$ & $0$  & No  \\
\ref{t010tr}.2 & $      $ & $      $ & $          $ &$ 0  $ & $0.5 $ & $1$  & No  \\
\ref{t010tr}.3 & $P_0(T)$ & $P_1(e)$ & $[P_{0}(T)]^2$ &$ 1 $   & $0$ & $0$  & No  \\
\ref{t010tr}.4 & $      $ & $      $ & $          $ &$\infty$ &$-\infty$ & $-\infty$  &No \\
\hline
\ref{t010tr}.5 & $      $ & $      $ & $          $ &$ -1 $& $1$ & $2$  &  No  \\
\ref{t010tr}.6 & $      $ & $      $ & $          $ &$ 0 $& $1.5$ & $2$   & No   \\
\ref{t010tr}.7 & $P_1(T)$ & $P_2(e)$ & $[P_{1}(T)]^2$ &$ 1 $& $  1$ & $1$   & No  \\
\ref{t010tr}.8 & $      $ & $      $ & $          $ &$\infty$& $-\infty$ & $-\infty$  &No \\
\hline
\ref{t010tr}.9 & $      $ & $      $ & $          $ &$ -1 $& $2$ & $3$  &  No  \\
\ref{t010tr}.10& $      $ & $      $ & $          $ &$ 0 $& $2.5$ & $3$   & No   \\
\ref{t010tr}.11& $P_2(T)$ & $P_3(e)$ & $[P_2(T)]^2$ &$ 1 $& $2$ & $2$   & No   \\
\ref{t010tr}.12& $      $ & $      $ & $          $ &$\infty$& $-\infty$ & $-\infty$  &No \\
\hline
\end{tabular}
\end{center}
\end{table}

Table \ref{t011tr}  demonstrates the convergence rates for  $(P_k(T),P_{k+1}(e),[P_{k+1}(T)]^2)$ with a stabilizer of different $j$ on triangular mesh.

\begin{table}[H]
 \caption{Element $(P_k(T),P_{k+1}(e),[P_{k+1}(T)]^2)$ on triangular mesh,
  $\3bar\cdot\3bar=O(h^{r_1})$ and $\|\cdot\|=O(h^{r_2})$.}
    \label{t011tr}
\begin{center}
\begin{tabular}{c|c|c|c|c|c|c|c}
\hline
element & $P_k(T)$ & $P_{k+1}(e)$ &  $[P_{k+1}(T)]^2$  & $j$ & $r_1$  & $r_2$ & Proved \\
\hline
\ref{t011tr}.1 & $      $ & $      $ & $          $ &$ -1 $ & $0$ & $0$  & No  \\
\ref{t011tr}.2 & $      $ & $      $ & $          $ &$ 0  $ & $1 $ & $1$  & No  \\
\ref{t011tr}.3 & $P_0(T)$ & $P_1(e)$ & $[P_{1}(T)]^2$ &$ 1 $   & $2$ & $2$  & No  \\
\ref{t011tr}.4 & $      $ & $      $ & $          $ &$\infty$ &$2$ & $2$  &Yes \\
\hline
\ref{t011tr}.5 & $      $ & $      $ & $          $ &$ -1 $& $1$ & $2$  &  No  \\
\ref{t011tr}.6 & $      $ & $      $ & $          $ &$ 0 $& $ 2$ & $3$   & No   \\
\ref{t011tr}.7 & $P_1(T)$ & $P_2(e)$ & $[P_{2}(T)]^2$ &$ 1 $& $  3$ & $4$   & No  \\
\ref{t011tr}.8 & $      $ & $      $ & $          $ &$\infty$& $3$ & $4$  &Yes \\
\hline
\ref{t011tr}.9 & $      $ & $      $ & $          $ &$ -1 $& $2$ & $3$  &  No  \\
\ref{t011tr}.10& $      $ & $      $ & $          $ &$ 0 $& $3$ & $4$   & No   \\
\ref{t011tr}.11& $P_2(T)$ & $P_3(e)$ & $[P_3(T)]^2$ &$ 1 $& $4$ & $5$   & No   \\
\ref{t011tr}.12& $      $ & $      $ & $          $ &$\infty$& $4$ & $5$  &Yes \\
\hline
\end{tabular}
\end{center}
\end{table}

\begin{remark}
The WG element $(P_k(T), P_{k+1}(e), [P_{k+1}(T)]^2)$ has order two supercloseness on triangular mesh.
\end{remark}

\smallskip

Table \ref{t012tr}  demonstrates the convergence rates for  $(P_k(T),P_{k+1}(e),[P_{k+2}(T)]^2)$ with a stabilizer of different $j$ on triangular mesh.

\begin{table}[H]
 \caption{Element $(P_k(T),P_{k+1}(e),[P_{k+2}(T)]^2)$ on triangular mesh,
  $\3bar\cdot\3bar=O(h^{r_1})$ and $\|\cdot\|=O(h^{r_2})$.}
    \label{t012tr}
\begin{center}
\begin{tabular}{c|c|c|c|c|c|c|c}
\hline
element & $P_k(T)$ & $P_{k+1}(e)$ &  $[P_{k+2}(T)]^2$  & $j$ & $r_1$  & $r_2$ & Proved \\
\hline
\ref{t012tr}.1 & $      $ & $      $ & $          $ &$ -1 $ & $0$ & $0$  & No  \\
\ref{t012tr}.2 & $      $ & $      $ & $          $ &$ 0  $ & $0 $ & $0$  & No  \\
\ref{t012tr}.3 & $P_0(T)$ & $P_1(e)$ & $[P_{2}(T)]^2$ &$ 1 $   & $0$ & $0$  & No  \\
\ref{t012tr}.4 & $      $ & $      $ & $          $ &$\infty$ &$0$ & $0$  &No \\
\hline
\ref{t012tr}.5 & $      $ & $      $ & $          $ &$ -1 $& $1$ & $2$  &  No  \\
\ref{t012tr}.6 & $      $ & $      $ & $          $ &$ 0 $& $ 1$ & $2$   & No   \\
\ref{t012tr}.7 & $P_1(T)$ & $P_2(e)$ & $[P_{3}(T)]^2$ &$ 1 $& $  1$ & $2$   & No  \\
\ref{t012tr}.8 & $      $ & $      $ & $          $ &$\infty$& $1$ & $2$  &No \\
\hline
\ref{t012tr}.9 & $      $ & $      $ & $          $ &$ -1 $& $2$ & $3$  &  No  \\
\ref{t012tr}.10& $      $ & $      $ & $          $ &$ 0 $& $2$ & $3$   & No   \\
\ref{t012tr}.11& $P_2(T)$ & $P_3(e)$ & $[P_4(T)]^2$ &$ 1 $& $2$ & $3$   & No   \\
\ref{t012tr}.12& $      $ & $      $ & $          $ &$\infty$& $2$ & $3$  &No \\
\hline
\end{tabular}
\end{center}
\end{table}

Table \ref{t01RT0tr}  demonstrates the convergence rates for  $(P_k(T),P_{k+1}(e),RT_k(T))$ with a stabilizer of different $j$ on triangular mesh.

\begin{table}[H]
 \caption{Element $(P_k(T),P_{k+1}(e),RT_k(T))$ on triangular mesh,
  $\3bar\cdot\3bar=O(h^{r_1})$ and $\|\cdot\|=O(h^{r_2})$.}
    \label{t01RT0tr}
\begin{center}
\begin{tabular}{c|c|c|c|c|c|c|c}
\hline
element & $P_k(T)$ & $P_{k+1}(e)$ &  $RT_k(T)$  & $j$ & $r_1$  & $r_2$ & Proved \\
\hline
\ref{t01RT0tr}.1 & $      $ & $      $ & $          $ &$ -1 $ & $0$ & $0$  & No  \\
\ref{t01RT0tr}.2 & $      $ & $      $ & $          $ &$ 0  $ & $0.5$ & $1$  & No  \\
\ref{t01RT0tr}.3 & $P_0(T)$ & $P_1(e)$ & $RT_0(T)$ &$ 1 $   & $1$ & $2$  & No  \\
\ref{t01RT0tr}.4 & $      $ & $      $ & $          $ &$\infty$ &$1$ & $2$  &No \\
\hline
\ref{t01RT0tr}.5 & $      $ & $      $ & $          $ &$ -1 $& $1$ & $2$  & No  \\
\ref{t01RT0tr}.6 & $      $ & $      $ & $          $ &$ 0 $& $ 1.5$ & $3$   & No   \\
\ref{t01RT0tr}.7 & $P_1(T)$ & $P_2(e)$ & $RT_1(T)$ &$ 1 $& $  2$ & $3$   & No  \\
\ref{t01RT0tr}.8 & $      $ & $      $ & $          $ &$\infty$&$-\infty$ & $-\infty$  &No\\
\hline
\ref{t01RT0tr}.9 & $      $ & $      $ & $          $ &$ -1 $& $2$ & $3$  & No  \\
\ref{t01RT0tr}.10& $      $ & $      $ & $          $ &$ 0 $& $2.5$ & $4$   & No   \\
\ref{t01RT0tr}.11& $P_2(T)$ & $P_3(e)$ & $RT_2(T)$ &$ 1 $& $3$ & $4$   & No   \\
\ref{t01RT0tr}.12& $      $ & $      $ & $          $ &$\infty$& $-\infty$ & $-\infty$  &No \\
\hline
\end{tabular}
\end{center}
\end{table}

Table \ref{t01RT1tr}  demonstrates the convergence rates for  $(P_k(T),P_{k+1}(e),RT_{k+1}(T))$ with a stabilizer of different $j$ on triangular mesh.

\begin{table}[H]
 \caption{Element $(P_k(T),P_{k+1}(e),RT_{k+1}(T))$ on triangular mesh,
  $\3bar\cdot\3bar=O(h^{r_1})$ and $\|\cdot\|=O(h^{r_2})$.}
    \label{t01RT1tr}
\begin{center}
\begin{tabular}{c|c|c|c|c|c|c|c}
\hline
element & $P_k(T)$ & $P_{k+1}(e)$ &  $RT_{k+1}(T)$  & $j$ & $r_1$  & $r_2$ & Proved \\
\hline
\ref{t01RT1tr}.1 & $      $ & $      $ & $          $ &$ -1 $ & $0$ & $0$  & No  \\
\ref{t01RT1tr}.2 & $      $ & $      $ & $          $ &$ 0  $ & $0$ & $0$  & No  \\
\ref{t01RT1tr}.3 & $P_0(T)$ & $P_1(e)$ & $RT_1(T)$ &$ 1 $   & $0$ & $0$  & No  \\
\ref{t01RT1tr}.4 & $      $ & $      $ & $          $ &$\infty$ &$0$ & $0$  &No \\
\hline
\ref{t01RT1tr}.5 & $      $ & $      $ & $          $ &$ -1 $& $1$ & $2$  & No  \\
\ref{t01RT1tr}.6 & $      $ & $      $ & $          $ &$ 0 $& $ 1$ & $2$   & No   \\
\ref{t01RT1tr}.7 & $P_1(T)$ & $P_2(e)$ & $RT_2(T)$ &$ 1 $& $  1$ & $2$   & No  \\
\ref{t01RT1tr}.8 & $      $ & $      $ & $          $ &$\infty$&$1$ & $2$  &No\\
\hline
\ref{t01RT1tr}.9 & $      $ & $      $ & $          $ &$ -1 $& $2$ & $3$  & No  \\
\ref{t01RT1tr}.10& $      $ & $      $ & $          $ &$ 0 $& $2$ & $3$   & No   \\
\ref{t01RT1tr}.11& $P_2(T)$ & $P_3(e)$ & $RT_3(T)$ &$ 1 $& $2$ & $3$   & No   \\
\ref{t01RT1tr}.12& $      $ & $      $ & $          $ &$\infty$& $2$ & $3$  &No \\
\hline
\end{tabular}
\end{center}
\end{table}

\section{The WG elements for $\ell>s$ on rectangular mesh}

Table \ref{t100rec}  demonstrates the convergence rates
    for  $(P_{k}(T),P_{k-1}(e),[P_{k-1}(T)]^2)$ with a stabilizer of different
     $j$ on rectangular mesh.

\begin{table}[H]
 \caption{Element $(P_{k}(T),P_{k-1}(e),[P_{k-1}(T)]^2)$ on rectangular mesh,
  $\3bar\cdot\3bar=O(h^{r_1})$ and $\|\cdot\|=O(h^{r_2})$.}
    \label{t100rec}
\begin{center}
\begin{tabular}{c|c|c|c|c|c|c|c}
\hline
element & $P_{k}(T)$ & $P_{k-1}(e)$ &  $[P_{k-1}(T)]^d$  & $j$ & $r_1$  & $r_2$ & Proved \\
\hline
\ref{t100rec}.1 & $      $ & $      $ & $          $ & $-1$ & $1$ & $2$ &Yes  \\
\ref{t100rec}.2 & $      $ & $      $ & $          $ & $0$  & $0.5$ & $1$  &No  \\
\ref{t100rec}.3 & $P_1(T)$ & $P_0(e)$ & $[P_0(T)]^2$ & $1$ & $0$ & $0$ &No  \\
\ref{t100rec}.4 & $      $ & $      $ & $          $ & $\infty$ & $-\infty$ & $-\infty$&No  \\
\hline
\ref{t100rec}.5 & $      $ & $      $ & $          $ & $-1$ & $2$ & $3$   &Yes   \\
\ref{t100rec}.6 & $      $ & $      $ & $          $ & 0 & $ 1.5 $ & $2 $  &No  \\
\ref{t100rec}.7 & $P_2(T)$ & $P_1(e)$ & $[P_1(T)]^2$ & $1$ & $ 1 $ & $1 $  &No  \\
\ref{t100rec}.8 & $      $ & $      $ & $          $ &  $\infty$ & $-\infty$ & $-\infty$&No  \\
\hline
\ref{t100rec}.9  & $      $ & $      $ & $          $ &$-1$ & $3$ & $4$ &Yes   \\
\ref{t100rec}.10 & $      $ & $      $ & $          $ & 0 & $2.5$ & $3$   &No  \\
\ref{t100rec}.11 & $P_3(T)$ & $P_2(e)$ & $[P_2(T)]^2$ & 1 & $2$ & $2$   &No \\
\ref{t100rec}.12 & $      $ & $      $ & $          $ & $\infty$ & $-\infty$ & $-\infty$&No  \\
\hline
\end{tabular}
\end{center}
\end{table}

Table \ref{t10RT0rec} demonstrates the convergence rates for
  $(P_{k}(T),P_{k-1}(e),[P_{k}(T)]^2)$ with a stabilizer of different $j$ on rectangular mesh.

\begin{table}[H]
 \caption{Element $(P_{k}(T),P_{k-1}(e),[P_{k}(T)]^2)$ on rectangular mesh,
  $\3bar\cdot\3bar=O(h^{r_1})$ and $\|\cdot\|=O(h^{r_2})$.}
    \label{t10RT0rec}
\begin{center}
\begin{tabular}{c|c|c|c|c|c|c|c}
\hline
element & $P_{k}(T)$ & $P_{k-1}(e)$ &  $[P_{k}(T)]^d$  & $j$ & $r_1$  & $r_2$ & Proved \\
\hline
\ref{t10RT0rec}.1 & $      $ & $      $ & $          $ & $-1$ & $1$ & $2$ & No  \\
\ref{t10RT0rec}.2 & $      $ & $      $ & $          $ & $0$  & $1.5$ & $2$  &No  \\
\ref{t10RT0rec}.3 & $P_1(T)$ & $P_0(e)$ & $[P_1(T)]^2$ & $1$ & $1$ & $1$ &No  \\
\ref{t10RT0rec}.4 & $      $ & $      $ & $          $ & $\infty$ & $-\infty$ & $-\infty$&No  \\
\hline
\ref{t10RT0rec}.5 & $      $ & $      $ & $          $ & $-1$ & $2$ & $3$   & No   \\
\ref{t10RT0rec}.6 & $      $ & $      $ & $          $ & 0 & $ 2.5 $ & $ 3 $  &No  \\
\ref{t10RT0rec}.7 & $P_2(T)$ & $P_1(e)$ & $[P_2(T)]^2$ & $1$ & $ 2 $ & $2 $  &No  \\
\ref{t10RT0rec}.8 & $      $ & $      $ & $          $ & $\infty$ & $-\infty$ & $-\infty$&No  \\
\hline
\ref{t10RT0rec}.9 & $      $ & $      $ & $          $ & $-1$ & $3$ & $4$ & No   \\
\ref{t10RT0rec}.10 & $      $ & $      $ & $          $ & 0 & $3.5$ & $4$   &No  \\
\ref{t10RT0rec}.11 & $P_3(T)$ & $P_2(e)$ & $[P_3(T)]^2$ & 1 & $3$ & $3$   &No \\
\ref{t10RT0rec}.12 & $      $ & $      $ & $          $ & $\infty$ & $-\infty$ & $-\infty$&No  \\
\hline
\end{tabular}
\end{center}
\end{table}

\begin{remark}
The WG element \ref{t10RT0rec}.1  achieves optimal convergence rates on triangular mesh while Theorem 4.9 in \cite{wwzz} predict only suboptimal convergence rate.
\end{remark}

\smallskip

Table \ref{t102rec}
 demonstrates the convergence rates for
  $(P_{k}(T),P_{k-1}(e),[P_{k+1}(T)]^2)$ with a stabilizer of different $j$ on rectangular mesh.

\begin{table}[H]
 \caption{Element $(P_{k}(T),P_{k-1}(e),[P_{k+1}(T)]^2)$ on rectangular mesh,
  $\3bar\cdot\3bar=O(h^{r_1})$ and $\|\cdot\|=O(h^{r_2})$.}
    \label{t102rec}
\begin{center}
\begin{tabular}{c|c|c|c|c|c|c|c}
\hline
element & $P_{k}(T)$ & $P_{k-1}(e)$ &  $[P_{k+1}(T)]^d$  & $j$ & $r_1$  & $r_2$ & Proved \\
\hline
\ref{t102rec}.1 & $      $ & $      $ & $          $ & $-1$ & $0$ & $0$ &  No  \\
\ref{t102rec}.2 & $      $ & $      $ & $          $ & $0$  & $0$ & $0$  &No  \\
\ref{t102rec}.3 & $P_1(T)$ & $P_0(e)$ & $[P_2(T)]^2$ & $1$ & $0$ & $0$ &No  \\
\ref{t102rec}.4 & $      $ & $      $ & $          $ & $\infty$ & $0$ & $0$&No   \\
\hline
\ref{t102rec}.5 & $      $ & $      $ & $          $ & $-1$ & $1$ & $2$   &  No   \\
\ref{t102rec}.6 & $      $ & $      $ & $          $ & 0 & $ 1 $ & $ 2 $  &No  \\
\ref{t102rec}.7 & $P_2(T)$ & $P_1(e)$ & $[P_3(T)]^2$ & $1$ & $ 1 $ & $2 $  &No  \\
\ref{t102rec}.8 & $      $ & $      $ & $          $ &  $\infty$ & $1$ & $2$&No   \\
\hline
\ref{t102rec}.9  & $      $ & $      $ & $          $ & $-1$ & $2$ & $3$ &  No   \\
\ref{t102rec}.10 & $      $ & $      $ & $          $ & 0 & $2$ & $3$   &No  \\
\ref{t102rec}.11 & $P_3(T)$ & $P_2(e)$ & $[P_4(T)]^2$ & 1 & $2$ & $3$   &No \\
\ref{t102rec}.12 & $      $ & $      $ & $          $ & $\infty$ & $2$ & $3$&No   \\
\hline
\end{tabular}
\end{center}
\end{table}

\begin{remark}
The numerical results in Table \ref{t102rec} show that the WG element $(P_{k}(T),P_{k-1}(e),[P_{k+1}(T)]^2)$ has suboptimal convergence rates on rectangular mesh. A new stabilizer free WG method is proposed  in \cite{Ye-Zhang-k-1} for the element $(P_{k}(T)$, $P_{k-1}(e)$,
  $[P_{k+1}(T)]^2)$ with optimal convergence rate, on general polygonal meshes.
\end{remark}

\smallskip

Table \ref{t10RTrec}
 demonstrates the convergence rates for
  $(P_{k}(T),P_{k-1}(e),RT_{k-1}(T))$ with a stabilizer of different $j$ on rectangular mesh.

\begin{table}[H]
 \caption{Element $(P_{k}(T),P_{k-1}(e),RT_{k-1}(T))$ on rectangular mesh,
  $\3bar\cdot\3bar=O(h^{r_1})$ and $\|\cdot\|=O(h^{r_2})$.}
    \label{t10RTrec}
\begin{center}
\begin{tabular}{c|c|c|c|c|c|c|c}
\hline
element & $P_{k}(T)$ & $P_{k-1}(e)$ &  $RT_{k-1}(T)$  & $j$ & $r_1$  & $r_2$ & Proved \\
\hline
\ref{t10RTrec}.1 & $      $ & $      $ & $          $ & $-1$ & $1$ & $2$ & No  \\
\ref{t10RTrec}.2 & $      $ & $      $ & $          $ & $0$  & $1.5$ & $2$  &No  \\ 
\ref{t10RTrec}.3 & $P_1(T)$ & $P_0(e)$ & $RT_{0}(T)$ & $1$ & $1$ & $1$ &No  \\
\ref{t10RTrec}.4 & $      $ & $      $ & $          $ & $\infty$ & $-\infty$ & $-\infty$&No  \\
\hline
\ref{t10RTrec}.5 & $      $ & $      $ & $          $ & $-1$ & $2$ & $3$   & No   \\
\ref{t10RTrec}.6 & $      $ & $      $ & $          $ & 0 &   $ 1.5 $ & $2 $  &No  \\ 
\ref{t10RTrec}.7 & $P_2(T)$ & $P_1(e)$ & $RT_{1}(T)$ & $1$ & $ 1 $ & $1 $  &No  \\
\ref{t10RTrec}.8 & $      $ & $      $ & $          $ &  $\infty$ & $-\infty$ & $-\infty$&No  \\
\hline
\ref{t10RTrec}.9 & $      $ & $      $ & $          $ & $-1$ & $3$ & $4$ &No   \\
\ref{t10RTrec}.10 & $      $ & $      $ & $          $ & 0 & $2.5$ & $3 $   &No  \\
\ref{t10RTrec}.11 & $P_3(T)$ & $P_2(e)$ & $RT_{2}(T)$ & 1 & $2$ & $2$   &No \\
\ref{t10RTrec}.12 & $      $ & $      $ & $          $ & $\infty$ & $-\infty$ & $-\infty$&No  \\
\hline
\end{tabular}
\end{center}
\end{table}

Table \ref{t10RT1rec}
 demonstrates the convergence rates for
  $(P_{k}(T),P_{k-1}(e),RT_{k}(T))$ with a stabilizer of different $j$ on rectangular mesh.

\begin{table}[H]
 \caption{Element $(P_{k}(T),P_{k-1}(e),RT_{k}(T))$ on rectangular mesh,
  $\3bar\cdot\3bar=O(h^{r_1})$ and $\|\cdot\|=O(h^{r_2})$.}
    \label{t10RT1rec}
\begin{center}
\begin{tabular}{c|c|c|c|c|c|c|c}
\hline
element & $P_{k}(T)$ & $P_{k-1}(e)$ &  $RT_{k}(T)$  & $j$ & $r_1$  & $r_2$ & Proved \\
\hline
\ref{t10RT1rec}.1 & $      $ & $      $ & $          $ & $-1$ & $0$ & $0$ & No  \\
\ref{t10RT1rec}.2 & $      $ & $      $ & $          $ & $0$  & $0$ & $1$  &No  \\
\ref{t10RT1rec}.3 & $P_1(T)$ & $P_0(e)$ & $RT_{1}(T)$ & $1$ & $0$ & $1$ &No  \\
\ref{t10RT1rec}.4 & $      $ & $      $ & $          $ & $\infty$ & $0$ & $1$&No  \\
\hline
\ref{t10RT1rec}.5 & $      $ & $      $ & $          $ & $-1$ & $1$ & $2$   & No   \\
\ref{t10RT1rec}.6 & $      $ & $      $ & $          $ & 0 &       $ 1 $ & $2 $  &No  \\
\ref{t10RT1rec}.7 & $P_2(T)$ & $P_1(e)$ & $RT_{2}(T)$ & $1$ &     $ 1 $ & $2 $  &No  \\
\ref{t10RT1rec}.8 & $      $ & $      $ & $          $ &  $\infty$ & $1$ & $2$&No  \\
\hline
\ref{t10RT1rec}.9  & $      $ & $      $ & $          $ & $-1$     & $2$ & $3$ &No   \\
\ref{t10RT1rec}.10 & $      $ & $      $ & $          $ & 0 &       $2$ & $3 $   &No  \\
\ref{t10RT1rec}.11 & $P_3(T)$ & $P_2(e)$ & $RT_{3}(T)$ & 1 &       $2$ & $3$   &No \\
\ref{t10RT1rec}.12 & $      $ & $      $ & $          $ & $\infty$ & $2$ & $3$&No  \\
\hline
\end{tabular}
\end{center}
\end{table}

\section{The WG elements for $\ell> s$ on triangular mesh}

Table \ref{t100tr}   demonstrates the convergence rates for  $(P_{k }(T),P_{k-1}(e),[P_{k-1}(T)]^2)$ with a stabilizer of different $j$ on triangular mesh.

\begin{table}[H]
 \caption{Element $(P_{k }(T),P_{k-1}(e),[P_{k-1}(T))$ on triangular mesh,
  $\3bar\cdot\3bar=O(h^{r_1})$ and $\|\cdot\|=O(h^{r_2})$.}
    \label{t100tr}
\begin{center}
\begin{tabular}{c|c|c|c|c|c|c|c}
\hline
element & $P_{k}(T)$ & $P_{k-1}(e)$ &  $[P_{k-1}(T)]^2$  & $j$ & $r_1$  & $r_2$ & Proved \\
\hline
\ref{t100tr}.1 & $      $ & $      $ & $          $ &$ -1 $ & $1$ & $2$  & Yes  \\
\ref{t100tr}.2 & $      $ & $      $ & $          $ &$ 0  $ & $0.5$ & $1$  & No  \\
\ref{t100tr}.3 & $P_1(T)$ & $P_0(e)$ & $[P_{0}(T)]^2$ &$ 1 $   & $0$ & $0$  & No  \\
\ref{t100tr}.4 & $      $ & $      $ & $          $ &$\infty$ &$-\infty$ & $-\infty$  &No \\
\hline
\ref{t100tr}.5 & $      $ & $      $ & $          $ &$ -1 $& $2$ & $3$  & Yes  \\
\ref{t100tr}.6 & $      $ & $      $ & $          $ &$ 0 $& $1.5$ & $2$   & No   \\
\ref{t100tr}.7 & $P_2(T)$ & $P_1(e)$ & $[P_{1}(T)]^2$ &$ 1 $& $  1$ & $1$   & No  \\
\ref{t100tr}.8 & $      $ & $      $ & $          $ &$\infty$& $-\infty$ & $-\infty$  &No \\
\hline
\ref{t100tr}.9 & $      $ & $      $ & $          $ &$ -1 $& $3$ & $4$  & Yes  \\
\ref{t100tr}.10& $      $ & $      $ & $          $ &$ 0 $& $2.5$ & $3$   & No   \\
\ref{t100tr}.11& $P_3(T)$ & $P_2(e)$ & $[P_2(T)]^2$ &$ 1 $& $2$ & $2$   & No   \\
\ref{t100tr}.12& $      $ & $      $ & $          $ &$\infty$& $-\infty$ & $-\infty$  &No \\
\hline
\end{tabular}
\end{center}
\end{table}

Table \ref{t101tr}   demonstrates the convergence rates for  $(P_{k }(T),P_{k-1}(e),[P_{k}(T)]^2)$ with a stabilizer of different $j$ on triangular mesh.

\begin{table}[H]
 \caption{Element $((P_{k}(T),P_{k-1}(e),[P_{k}(T)]^2))$ on triangular mesh,
  $\3bar\cdot\3bar=O(h^{r_1})$ and $\|\cdot\|=O(h^{r_2})$.}
    \label{t101tr}
\begin{center}
\begin{tabular}{c|c|c|c|c|c|c|c}
\hline
element & $P_{k}(T)$ & $P_{k-1}(e)$ &  $[P_{k}(T)]^2$  & $j$ & $r_1$  & $r_2$ & Proved \\
\hline
\ref{t101tr}.1 & $      $ & $      $ & $          $ &$ -1 $ & $0$ & $0$  &  No  \\
\ref{t101tr}.2 & $      $ & $      $ & $          $ &$ 0  $ & $0 $ & $0 $  & No  \\
\ref{t101tr}.3 & $P_1(T)$ & $P_0(e)$ & $[P_{1}(T)]^2$ &$ 1 $   & $0$ & $0$  & No  \\
\ref{t101tr}.4 & $      $ & $      $ & $          $ &$\infty$ &$-\infty$ & $-\infty$  &No \\
\hline
\ref{t101tr}.5 & $      $ & $      $ & $          $ &$ -1 $& $1$ & $2$  &  No  \\
\ref{t101tr}.6 & $      $ & $      $ & $          $ &$ 0 $& $1 $ & $2$   & No   \\
\ref{t101tr}.7 & $P_2(T)$ & $P_1(e)$ & $[P_{2}(T)]^2$ &$ 1 $& $  1$ & $2$   & No  \\
\ref{t101tr}.8 & $      $ & $      $ & $          $ &$\infty$& $-\infty$ & $-\infty$  &No \\
\hline
\ref{t101tr}.9 & $      $ & $      $ & $          $ &$ -1 $& $2$ & $3$  &  No  \\
\ref{t101tr}.10& $      $ & $      $ & $          $ &$ 0 $& $2$ & $3$   & No   \\
\ref{t101tr}.11& $P_3(T)$ & $P_2(e)$ & $[P_3(T)]^2$ &$ 1 $& $2$ & $3$   & No   \\
\ref{t101tr}.12& $      $ & $      $ & $          $ &$\infty$& $-\infty$ & $-\infty$  &No \\
\hline
\end{tabular}
\end{center}
\end{table}

\begin{remark}
The WG element $(P_k(T), P_{k-1}(e), [P_{k}(T)]^2)$ performs  better on rectangular mesh than on triangular mesh.
\end{remark}

\smallskip

Table \ref{t102tr}   demonstrates the convergence rates for  $(P_{k}(T),P_{k-1}(e),[P_{k+1}(T)]^2)$ with a stabilizer of different $j$ on triangular mesh.

\begin{table}[H]
 \caption{Element $((P_{k }(T),P_{k-1}(e),[P_{k+1}(T)]^2))$ on triangular mesh,
  $\3bar\cdot\3bar=O(h^{r_1})$ and $\|\cdot\|=O(h^{r_2})$.}
    \label{t102tr}
\begin{center}
\begin{tabular}{c|c|c|c|c|c|c|c}
\hline
element & $P_{k }(T)$ & $P_{k-1}(e)$ &  $[P_{k+1}(T)]^2$  & $j$ & $r_1$  & $r_2$ & Proved \\
\hline
\ref{t102tr}.1 & $      $ & $      $ & $          $ &$ -1 $ & $0$ & $0$  &  No  \\
\ref{t102tr}.2 & $      $ & $      $ & $          $ &$ 0  $ & $0 $ & $0 $  & No  \\
\ref{t102tr}.3 & $P_1(T)$ & $P_0(e)$ & $[P_{2}(T)]^2$ &$ 1 $   & $0$ & $0$  & No  \\
\ref{t102tr}.4 & $      $ & $      $ & $          $ &$\infty$ &$0$ & $0$  &No  \\
\hline
\ref{t102tr}.5 & $      $ & $      $ & $          $ &$ -1 $& $1$ & $2$  &  No  \\
\ref{t102tr}.6 & $      $ & $      $ & $          $ &$ 0 $& $1 $ & $2$   & No   \\
\ref{t102tr}.7 & $P_2(T)$ & $P_1(e)$ & $[P_{3}(T)]^2$ &$ 1 $& $  1$ & $2$   & No \\
\ref{t102tr}.8 & $      $ & $      $ & $          $ &$\infty$& $1$ & $2$  &No  \\
\hline
\ref{t102tr}.9 & $      $ & $      $ & $          $ &$ -1 $& $2$ & $3$  &  No  \\
\ref{t102tr}.10& $      $ & $      $ & $          $ &$ 0 $& $2$ & $3$   & No   \\
\ref{t102tr}.11& $P_3(T)$ & $P_2(e)$ & $[P_4(T)]^2$ &$ 1 $& $2$ & $3$   & No   \\
\ref{t102tr}.12& $      $ & $      $ & $          $ &$\infty$& $2$ & $3$  &No  \\
\hline
\end{tabular}
\end{center}
\end{table}

Table \ref{t10RT0tr}   demonstrates the convergence rates for  $(P_{k }(T),P_{k-1}(e),RT_{k-1}(T))$ with a stabilizer of different $j$ on triangular mesh.

\begin{table}[H]
 \caption{Element $( P_{k }(T),P_{k-1}(e),RT_{k-1}(T))$ on triangular mesh,
  $\3bar\cdot\3bar=O(h^{r_1})$ and $\|\cdot\|=O(h^{r_2})$.}
    \label{t10RT0tr}
\begin{center}
\begin{tabular}{c|c|c|c|c|c|c|c}
\hline
element & $P_{k}(T)$ & $P_{k-1}(e)$ &  $RT_{k-1}(T)$  & $j$ & $r_1$  & $r_2$ & Proved \\
\hline
\ref{t10RT0tr}.1 & $      $ & $      $ & $          $ &$ -1 $ & $1$ & $2$  & No  \\
\ref{t10RT0tr}.2 & $      $ & $      $ & $          $ &$ 0  $ & $1$ & $2$  & No  \\
\ref{t10RT0tr}.3 & $P_1(T)$ & $P_0(e)$ & $RT_0(T)$ &$ 1 $   & $1$ & $1$  & No  \\
\ref{t10RT0tr}.4 & $      $ & $      $ & $          $ &$\infty$ &$-\infty$ & $-\infty$  &No \\
\hline
\ref{t10RT0tr}.5 & $      $ & $      $ & $          $ &$ -1 $& $2$ & $3$  & No  \\
\ref{t10RT0tr}.6 & $      $ & $      $ & $          $ &$ 0 $& $2 $ & $3$   & No   \\
\ref{t10RT0tr}.7 & $P_2(T)$ & $P_1(e)$ & $RT_1(T)$ &$ 1 $& $ 2$ & $2$   & No  \\
\ref{t10RT0tr}.8 & $      $ & $      $ & $          $ &$\infty$& $-\infty$ & $-\infty$  &No \\
\hline
\ref{t10RT0tr}.9 & $      $ & $      $ & $          $ &$ -1 $& $3$ & $4$  & No  \\
\ref{t10RT0tr}.10& $      $ & $      $ & $          $ &$ 0 $& $3$ & $4$   & No   \\
\ref{t10RT0tr}.11& $P_3(T)$ & $P_2(e)$ & $RT_2(T)$ &$ 1 $& $3$ & $3$   & No   \\
\ref{t10RT0tr}.12& $      $ & $      $ & $          $ &$\infty$& $-\infty$ & $-\infty$  &No \\
\hline
\end{tabular}
\end{center}
\end{table}

Table \ref{t10RT1tr}   demonstrates the convergence rates for  $(P_{k }(T),P_{k-1}(e),RT_{k}(T))$ with a stabilizer of different $j$ on triangular mesh.

\begin{table}[H]
 \caption{Element $(P_{k }(T),P_{k-1}(e),RT_{k}(T))$ on triangular mesh,
  $\3bar\cdot\3bar=O(h^{r_1})$ and $\|\cdot\|=O(h^{r_2})$.}
    \label{t10RT1tr}
\begin{center}
\begin{tabular}{c|c|c|c|c|c|c|c}
\hline
element & $P_{k}(T)$ & $P_{k-1}(e)$ &  $RT_{k}(T)$  & $j$ & $r_1$  & $r_2$ & Proved \\
\hline
\ref{t10RT1tr}.1 & $      $ & $      $ & $          $ &$ -1 $ & $0$ & $0$  & No  \\
\ref{t10RT1tr}.2 & $      $ & $      $ & $          $ &$ 0  $ & $0$ & $0$  & No  \\
\ref{t10RT1tr}.3 & $P_1(T)$ & $P_0(e)$ & $RT_1(T)$ &$ 1 $   & $0$ & $0$  & No  \\
\ref{t10RT1tr}.4 & $      $ & $      $ & $          $ &$\infty$ &$0$ & $0$  &No \\
\hline
\ref{t10RT1tr}.5 & $      $ & $      $ & $          $ &$ -1 $& $1$ & $2$  & No  \\
\ref{t10RT1tr}.6 & $      $ & $      $ & $          $ &$ 0 $& $1 $ & $2$   & No   \\
\ref{t10RT1tr}.7 & $P_2(T)$ & $P_1(e)$ & $RT_2(T)$ &$ 1 $& $ 1$ & $2$   & No  \\
\ref{t10RT1tr}.8 & $      $ & $      $ & $          $ &$\infty$& $1$ & $2$  &No \\
\hline
\ref{t10RT1tr}.9 & $      $ & $      $ & $          $ &$ -1 $& $2$ & $3$  & No  \\
\ref{t10RT1tr}.10& $      $ & $      $ & $          $ &$ 0 $& $2$ & $3$   & No   \\
\ref{t10RT1tr}.11& $P_3(T)$ & $P_2(e)$ & $RT_3(T)$ &$ 1 $& $2$ & $3$   & No   \\
\ref{t10RT1tr}.12& $      $ & $      $ & $          $ &$\infty$& $2$ & $3$  &No \\
\hline
\end{tabular}
\end{center}
\end{table}

\end{document}